\documentclass[10 pt, conference]{ieeeconf}
\IEEEoverridecommandlockouts                              

\usepackage{cite}
\expandafter\let\csname proof\endcsname\relax
\expandafter\let\csname endproof\endcsname\relax
\usepackage{amsmath,amssymb,amsfonts}
\usepackage{graphicx}
\usepackage{textcomp}
\usepackage{makecell}
\usepackage{xcolor}
\usepackage{amsthm}
\usepackage{hyperref}
\usepackage[ruled]{algorithm2e} 
\usepackage{algorithmic}
\usepackage{subfigure}
\usepackage{booktabs}
\usepackage{balance}
\def\BibTeX{{\rm B\kern-.05em{\sc i\kern-.025em b}\kern-.08em
    T\kern-.1667em\lower.7ex\hbox{E}\kern-.125emX}}
    \newtheorem{theorem}{Theorem}

\newtheorem{lemma}{Lemma}

\newtheorem*{remark}{Remark}
\newtheorem{definition}{Definition}

\usepackage{threeparttable}

\title{\LARGE \bf
Explicit Solutions for Safety Problems Using Control Barrier Functions
\vspace{-10pt}
}
\author{Han Wang, Kostas Margellos, Antonis Papachristodoulou
\thanks{The authors are with the Department of Engineering Science, University of Oxford, Oxford, United Kingdom. E-mails: {\tt\small \{han.wang, kostas.margellos, antonis\}@eng.ox.ac.uk}}
}

\begin{document}

\maketitle
\thispagestyle{empty}
\pagestyle{empty}
\begin{abstract}
The control Barrier function approach has been widely used for safe controller synthesis. By solving an online convex quadratic programming problem, an optimal safe controller can be synthesized implicitly in state-space. Since the solution is unique, the mapping from state-space to control inputs is injective, thus enabling us to evaluate the underlying relationship. In this paper we aim at explicitly synthesizing a safe control law as a function of the state for nonlinear control-affine systems with limited control ability. We propose to transform the online quadratic programming problem into an offline parameterized optimisation problem which considers states as parameters. The obtained explicit safe controller is shown to be a piece-wise Lipschitz continuous function over the partitioned state space if the program is feasible. We address the infeasible cases by solving a parameterized adaptive control Barrier function-based quadratic programming problem. Extensive simulation results show the state-space partition and the controller properties.
\end{abstract}

\section{Introduction}\label{sec:introduction}
Safety verification and safe controller design for dynamical systems have attracted significant attention in safety critical applications such as collision avoidance, traffic flow control, adaptive cruise control etc. Safety verification aims at verifying that the system trajectory belongs to a safe set over an infinite time horizon. Motivated by invariance analysis \cite{blanchini1999set}, safety is proven to be equivalent to set invariance \cite{prajna2005necessity}. To establish a relationship between invariant and safe sets, the Barrier certificates method has been proposed \cite{prajna2004safety}, \cite{prajna2007framework}. This formulation is shown to be efficient for both deterministic and stochastic dynamical system settings; safe controller design remains a challenge, especially for vector fields for which safety cannot be guaranteed without a control input. To address this issue, a control Barrier function approach was proposed \cite{ames2016control}.

The Control Barrier Function approach is motivated from the control Lyapunov function approach \cite{freeman1996control}, which guarantees stability of a system by imposing Lyapunov's conditions. Given system dynamics and a predefined Barrier function, safety is guaranteed by redirecting the vector field into the safe set. With an additional relaxation term encoding the distance from the boundary, the safe controller is guaranteed to be Lipschitz continuous \cite{xu2015robustness}. To trade-off safety and performance, a quadratic programming (QP) formulation is used to integrate control Barrier function and control Lyapunov function constraints \cite{ames2016control}. This approach was later extended to high-order cases \cite{xiao2019control}, \cite{tan2021high}, and used in many applications \cite{xu2017correctness}, \cite{chen2017obstacle}. 

The success of this methodology relies on a closed-loop safe controller obtained by solving the online point-wise QP. However, this framework also makes the performance of the controller hard to evaluate, since the closed loop response is implicit. Moreover, although solving a convex quadratic programming problem is computationally cheap, a QP solver is sometimes overloaded for embedded devices. This motivates us to explicitly synthesize a closed-loop safe controller, and characterize its properties.

By considering the states as parameters in the quadratic programming formulation, the optimisation problem turns out to be a multi-parametric one \cite{gal1972multiparametric}. Major advances have been made in this field for applications related to our problem, such as explicit MPC \cite{bemporad2002explicit}, \cite{tondel2003algorithm}. To solve this problem, the parameter space is decomposed into critical regions \cite{fiacco1990sensitivity}. Inside every critical region, the problem degenerates into an equality-constrained programming problem, whose explicit solution can be evaluated. For the class of quadratic programming and quadratically constrained quadratic programming problems, the critical regions can be enumerated in polynomial time \cite{pappas2021exact}. 

In this paper we consider the explicit controller design problem for a general nonlinear control-affine dynamical system with limited control ability. We parameterize the state-space, and reformulate the control Barrier function-based quadratic programming problem as a multi-parametric programming problem. By analyzing the sensitivity of the optimisation problem at given parameters, the state-space is partitioned into multiple critical regions. The safe control law is evaluated explicitly inside every critical region. We further consider the case where the original problem is infeasible for some states. In this case, an alternative adaptive control Barrier function formulation is used. Following a similar analysis procedure, we obtain both the safe controller and adapted relaxation term designs as explicit piece-wise functions. Previous work to analyse the explicit closed loop response \cite{jankovic2018robust}, \cite{tan2021undesired} concentrated on cases without control limitations and did not consider the analysis of Lipschitz continuity and feasibility.

The remainder of the paper is organized as follows. Preliminaries including control Barrier functions and sensitivity analysis are introduced in Section \ref{sec:preliminaries}. The explicit controller design approach is introduced in Section \ref{sec:cbf}. Numerical simulation results are reported in Section \ref{sec:simulation}. Section \ref{sec:conclusion} concludes the paper.

\section{Preliminaries and Notations}\label{sec:preliminaries}
\textit{Notation}: For matrix $A$, we use $A_i$ to represent the $i$-th row of $A$. For an index set $\mathcal{I}$, $A_\mathcal{I}$ denotes the matrix composed of vectors $A_i,i\in\mathcal{I}$. We use $A \succeq 0$ and $A \succ 0$ to denote that $A$ is positive semi-definite 
and positive definite, respectively. $\mathbb{R}$ denotes the real space, $\mathbb{R}_+$ is the positive half space, and $\mathbb{R}_+^m$ lifts the dimension to $m$.

\subsection{Control Barrier Functions}
Consider a nonlinear control-affine system
\begin{equation}\label{eq:nonlnsys}
    \dot x = f(x)+g(x)u,
\end{equation}
with $x \in\mathbb{R}^n$, $u \in\mathcal{U}\subset \mathbb{R}^m$, $f(x):\mathbb{R}^n\to\mathbb{R}^n,$ and $g(x):\mathbb{R}^n\to \mathbb{R}^{n\times m}$. Both functions are further assumed to be Lipschitz continuous. Defining the solution of \eqref{eq:nonlnsys} to be $\psi(u,t,x_0)$, where $x_0$ represents the initial condition and $t$ denotes time, our goal is to design a closed loop controller $u(x)$ so that $\psi(u(x),t,x_0)$ stays within a given safe set $\mathcal{S}$ for any $t$. The existence and uniqueness of the solutions $\psi(u(x),t,x_0)$ are guaranteed by assuming the system is \textit{forward complete}, i.e. $\psi(u(x),t,x_0)$ is unique  for any $t\ge 0$.  

The safe set $\mathcal{S}$ is represented by the zero-super level set of a continuously differentiable function $s(x)$. Dually, the unsafe set $\bar{\mathcal{S}}$ can be defined as the complementary set. For ease of illustration we give the algebraic expressions of safe set $\mathcal{S}$, boundary of the safe set $\partial{\mathcal{S}}$, interior of the safe set $\mathrm{Int}(\mathcal{S})$ and unsafe set $\bar{\mathcal{S}}$:
\begin{subequations}\label{eq:safe}
\begin{align}
    \mathcal{S}&:=\{x\in\mathbb{R}^n|s(x)\ge0\},\\
    \partial{\mathcal{S}}&:=\{x\in\mathbb{R}^n|s(x)=0\},\\
    \mathrm{Int}(\mathcal{S})&:=\{x\in\mathbb{R}^n|s(x)>0\},\\
    \bar{\mathcal{S}}&:=\{x\in\mathbb{R}^n|s(x)<0\}.
\end{align}
\end{subequations}

With this formulation, the safe control design problem boils down to finding $u(x)\in\mathcal{U}$, such that $s(\psi(u(x),t,x_0))\ge 0$ for any $t$. To achieve this, a control Barrier function-based quadratic programming approach was proposed \cite{ames2016control} .

Control Barrier functions are an extension to Barrier certificates \cite{prajna2004safety} for safety verification. It has been revealed in these papers that safety is closely related to the notion of \textit{forward invariance}.
\begin{definition}\label{def:invariance}
A set $\mathcal{B} \subset \mathbb{R}^n$ is said to be forward invariant with respect to vector field \eqref{eq:nonlnsys}, if for any $x_0\in\mathcal{B}$, there exists $u\in\mathcal{U}$ such that $\psi(u,t,x_0)\in\mathcal{B}$.
\end{definition}
The relationship between safety and forward invariance is demonstrated in the following equivalence lemma.
\begin{lemma}[\cite{liu2021converse}]\label{lem:equivalence}
System \eqref{eq:nonlnsys} is safe under $\mathcal{S}$, if and only if there exists a forward invariant set $\mathcal{B}\subseteq\mathcal{S}$.
\end{lemma}
Clearly, given a forward invariant set $\mathcal{B}$, a safe controller $u(x)$ always exists for any $x\in\mathcal{B}$. The control Barrier function approach answers the question of how to design a closed loop safe controller $u(x)$ inside $\mathcal{B}$, and how to guarantee the resulting safe controller is Lipschitz continuous. The notion of control Barrier functions is facilitated by the notion of extended class-$\mathcal{K}$ functions.
\begin{definition}\label{def:classk}
A continuous function $\alpha(\cdot):(-b,a)\to (-\infty,+\infty)$ is said to be an extended class-$\mathcal{K}$ function for positive $a$ and $b$, if it is strictly increasing and $\alpha(0)=0$. 
\end{definition}

\begin{definition}
For the control-affine dynamical system \eqref{eq:nonlnsys}, a continuously differentiable function $B(\cdot):\mathbb{R}^n\to \mathbb{R}$ is said to be a control Barrier function for the set $\mathcal{B}$, if there exists an extended class-$\mathcal{K}$ function $\alpha(\cdot)$ and a set $\mathcal{C}$, where $\mathcal{B}\subseteq\mathcal{C}\subset \mathbb{R}^n$, such that for any $x\in\mathcal{C}$,
\begin{equation}\label{eq:cbf}
   \mathop{\sup\limits_{u\in\mathcal{U}}}[\mathcal{L}_fB(x)+\mathcal{L}_gB(x)u+\alpha(B(x))]\ge0.
\end{equation}
Here $\mathcal{L}_fB(x)$ and $\mathcal{L}_gB(x)$ are Lie derivatives, which are defined by $\mathcal{L}_fB(x):=\frac{\partial B(x)}{\partial x}f(x)$ and $\mathcal{L}_gB(x):=\frac{\partial B(x)}{\partial x}g(x)$, respectively.
\end{definition}
Given a control Barrier function $B(x)$, the control admissible set corresponding to \eqref{eq:cbf} is defined by
\begin{equation}\label{eq:cbfset}
    K_{cbf}(x):=\{u\in\mathcal{U}:\mathcal{L}_fB(x)+\mathcal{L}_gB(x)u+\alpha(B(x))\ge0\}.
\end{equation}
\begin{theorem}{\cite{ames2016control}}
Consider a control Barrier function $B(x)$ defined on $\mathcal{C}$. Then for any $x\in\mathcal{C}$, any $u(x)\in K_{cbf}(x)$ will render the set $\mathcal{B}$ forward invariant.
\end{theorem}
Clearly, the set $K_{cbf}(x)$ maps state $x$ to control input $u$. This motivates us to compute the explicit expression of closed-loop safe controller design, possibly compromising performance.

\subsection{Sensitivity Analysis}
Consider the following multi-parametric optimisation problem
\begin{equation}\label{eq:mpop}
    \begin{split}
        \min_x~ &f(x,\theta)\\
        \mathrm{subject~to}~&g_i(x,\theta)\le 0,\forall i\in\mathbb{I},\\
        &h_j(x,\theta)=0,\forall j\in\mathbb{J},\\
    \end{split}
\end{equation}
where $\theta\in\mathbb{R}^v$ is a parameter vector. Suppose the functions $f(x,\theta),g_i(x,\theta),h_j(x,\theta),\forall i\in\mathbb{I},\forall j\in\mathbb{J}$ are convex in $x$. The problem here is to represent the optimizer of \eqref{eq:mpop} as a function of $\theta$, i.e. $x(\theta)$. To do this, one strategy is to solve the multi-parametric optimisation problem with different parameters $\theta^*$. Then the problem turns out to be a general optimisation problem. In a small region around $\theta^*$, $x(\theta)$ keeps a unifying expression. The readers are referred to \cite{fiacco1990sensitivity} for more information on this topic.

For a quadratic programming problem for which $f(x,\theta)$ is quadratic in $x$ and $h_j(x,\theta),\forall j\in\mathbb{J}$,  $g_i(x,\theta),\forall i\in\mathbb{I}$ are linear in $x$, for any $\theta^*$, there always exists a polyhedral neighborhood in which the expression $\eta(\theta)$ is unchanged \cite{blanchini1999set}. This motivates us to consider using polyhedral regions to partition the state space for explicit safe controller design.

\section{Explicit Safe Controller Design}\label{sec:cbf}
In this section, we show how to synthesize an explicit safe controller by evaluating the explicit solution to the control Barrier function-based quadratic programming problems. Suppose that for system \eqref{eq:nonlnsys} and a set $\mathcal{C}$, there exists a control Barrier function $B(x)$, with $\mathcal{B}$ defined by its zero-super level set.

\subsection{Control Barrier Functions Based QP}\label{sec:cbfqp}
To design a safe controller satisfying the control Barrier functions constraint \eqref{eq:cbf}, a QP is used to search for feasible points and to guarantee controller performance. Usually, a nominal controller $u^{\mathrm{des}}(x)\in\mathcal{U}$ is synthesized with other methods in advance, e.g. optimal control or PID, for controller performance. The goal for the quadratic programming setting is to find a safe controller $u^*(x)$ with a minimized bias from $u^{\mathrm{des}}(x)$. The QP is given in the following form:
\begin{equation}\label{eq:CBFQP}
\begin{split}
    u^*(x)&=\arg\min_{u}~\frac{1}{2}||u-u^{\mathrm{des}}(x)||_2^2\\
    &\mathrm{subject~to}~\mathcal{L}_fB(x)+\mathcal{L}_gB(x)u+\alpha(B(x))\ge 0,\\
    &~~~~~Au+b\le 0,
\end{split}
\end{equation}
where $A\in\mathbb{R}^{p\times m}$ and $b\in\mathbb{R}^p$. Here $Au+b\le 0$ represents the control input constraints. The Lagrangian of problem \eqref{eq:CBFQP} is
\begin{equation}\label{eq:cbflagrangian}
\begin{split}
 L(u,\lambda,\mu) = &\frac{1}{2}||u-u^{\mathrm{des}}||_2^2+\mu^\top (Au+b)\\
 -&\lambda(\mathcal{L}_fB(x)+\mathcal{L}_gB(x)u+\alpha(B(x))),
\end{split}
\end{equation}
where $\lambda\in\mathbb{R}_+,\mu\in\mathbb{R}_+^p$ are the Lagrange multipliers.

Clearly, in \eqref{eq:cbf} the objective function is strongly convex, and the constraints are linear, thus strong duality holds. The optimal controller $u^*(x)$ can be derived by applying the Karush-Kuhn-Tucker (KKT) conditions
\begin{subequations}\label{eq:cbfkkt}
\begin{align}
    &u^*(x)-u^{\mathrm{des}}+G(x)^\top \lambda +A^\top \mu=0,\label{eq:cbfkkt1}\\
    &\mu^\top(Au^*(x)+b)=0,\label{eq:cbfkkt2}\\
    &\lambda(F(x)+G(x)u^*(x)+\Lambda(x))=0,\label{eq:cbfkkt3}\\
    &\lambda\ge0,F(x)+G(x)u^*(x)+\Lambda(x)\le0,\label{eq:cbfkkt4}\\
    &\mu\ge0,Au^*(x)+b\le0.\label{eq:cbfkkt5}
\end{align}
\end{subequations}
Here we substitute $-\mathcal{L}_fB(x)$, $-\mathcal{L}_gB(x)$, $-\alpha(B(x))$ for $F(x)$, $G(x)$, $\Lambda(x)$, respectively. Though $\lambda$, $\mu$ are also functions of parameter $x$, we omit the argument without ambiguity. The reason why we are interested in the KKT condition is that, at a given state, the optimal control law $u^*(x)$ fulfills the equality conditions in \eqref{eq:cbfkkt1} -- \eqref{eq:cbfkkt3}, which enables us to explore explicit expressions. On the other hand, the inequality conditions \eqref{eq:cbfkkt4} -- \eqref{eq:cbfkkt5} specify optimality and feasibility conditions, which define critical regions in the sequel.

To ease demonstration, we assume that the relative degree of the control Barrier function $B(x)$ is one. The results of this paper can be extended to high order cases, and this is left for future work. Moreover, the primal degenerate cases where linearly dependent constraints exist, and dual degenerate cases where $\lambda=0$ or $\mu=0$, are omitted. These cases are addressed in Section \ref{sec:partition} via a pruning scheme.

For a given $x$, we separate the control admissible constraints into active constraints $A_iu^*(x)+b_i=0$, and inactive ones $A_ju^*(x)+b_j<0$. The active index set $\mathcal{I}(x)$ is defined by $\mathcal{I}(x)=\{i|A_iu^*(x)+b_i=0\}$, while the complement set of inactive indices is denoted by $\bar{\mathcal{I}}(x)$. At a given $x$, the index set can be verified by testing all the control admissible constraints. With such a separation scheme, we are able to evaluate the solution of \eqref{eq:mpop} from three cases by considering whether the control Barrier function constraint is active or not.

\textit{Case 1}: Given $x$, the control Barrier function constraint is inactive, i.e. $F(x)+G(x)u+\Lambda(x)<0.$ The optimizer of \eqref{eq:mpop} in this case is given by
\begin{equation}\label{eq:case1}
    u^*(x)=u^\mathrm{des}.
\end{equation}
To see this, we have already assumed that $u^\mathrm{des}\in\mathcal{U}$. The saddle point of the primal function $\frac{1}{2}||u-u^{\mathrm{des}}||_2^2$ is $u^{\mathrm{des}}$, which is a feasible point satisfying the constraint $Au+b\ge0$. Hence, the minimizer of problem \eqref{eq:mpop} when the control Barrier function constraint is inactive, is $u^\mathrm{des}$.

\textit{Case 2}: Given $x$, the control Barrier function constraint is active, i.e. $F(x)+G(x)u(x)+\Lambda(x)=0$. The control admissible constraints are all inactive i.e. $\mathcal{I}(x)=\emptyset$ and $\mu=0$. Suppose that the optimisation problem \eqref{eq:mpop} is feasible, then the optimizer $u^*(x)$ and multiplier $\lambda(x)$ in this case are
\begin{subequations}\label{eq:case2}
\begin{align}
    u^*(x) &= u^\mathrm{des} - \frac{G(x)^\top
    (F(x)+G(x)u^{\mathrm{des}}+\Lambda(x))}{G(x)G(x)^\top},\label{eq:case2u}\\
    \lambda(x) &= \frac{F(x)+G(x)u^{\mathrm{des}}+\Lambda(x)}{G(x)G(x)^\top}.\label{eq:case2lambda}
\end{align}
\end{subequations}

To see this, by substituting $\mu=0$ into condition 
\eqref{eq:cbfkkt1} we have:
\begin{equation}\label{eq:lem31}
    u^*(x) = u^\mathrm{des}-G(x)^\top \lambda(x).
\end{equation}
Here $\lambda(x) \ne 0$ since the control Barrier function constraint is active. Substitute $u^*(x)$ from \eqref{eq:lem31} into $F(x)+G(x)u^*(x)+\Lambda(x)=0$ to get \eqref{eq:case2lambda}, 
where $G(x)G(x)^\top \ne 0$ provided that the relative degree of the control Barrier function is one. We can then substitute $\lambda(x)$ from \eqref{eq:case2lambda} into \eqref{eq:lem31} to obtain \eqref{eq:case2u}.

\textit{Case 3}: Given $x$, the control Barrier function constraint is active, i.e. $F(x)+G(x)u+\Lambda(x)= 0,\lambda(x)>0$, the control limit constraints are active with the index set $\mathcal{I}(x)$. Suppose that the problem \eqref{eq:mpop} is feasible, then the optimizer and multipliers of \eqref{eq:cbflagrangian} in this case are
\begin{subequations}\label{eq:case3}
\begin{align}
    u^*(x) &= u^{\mathrm{des}}-G(x)^\top\lambda(x)-{A_{\mathcal{I}}}^\top \mu(x),\label{eq:case3eq1}\\
    \lambda(x) &= (G(x)G(x)^\top)^{-1}\label{eq:case3eq2}\\
    &\times(F(x)+G(x)u^{\mathrm{des}}-G(x)A_{\mathcal{I}}\mu+\Lambda(x)),\nonumber\\
    \mu(x)&=(A_{\mathcal{I}}\tilde G^+(x){A_{\mathcal{I}}}^\top)^{-1}\\
    \times&(A_{\mathcal{I}}\tilde G^-(x)u^{\mathrm{des}}-A_{\mathcal{I}}G(x)^\top(G(x)G(x)^\top)^{-1}F(x)-\nonumber\\
    &A_{\mathcal{I}}G(x)^{\top}(G(x)G(x)^\top)^{-1}\Lambda(x)+b_{\mathcal{I}}),\nonumber
\end{align}
\end{subequations}
where $\tilde G^-(x):= I-G(x)^\top(G(x)G(x)^\top)^{-1}G(x)$, $\tilde G^+(x):=I+G(x)^\top(G(x)G(x)^\top)^{-1}G(x)$.

The results are computed by applying the KKT condition \eqref{eq:cbfkkt1} and the active constraint equations at the saddle point $(u^*(x),\lambda(x),\mu(x))$; we have $u^*(x)-u^{\mathrm{des}}+G(x)^\top\lambda(x)+{{A_{\mathcal{I}}}^\top\mu(x)=0}$, $F(x)+G(x)u^*(x)+\Lambda(x)=0$, $A_{\mathcal{I}}u^*(x)+b_{\mathcal{I}}=0$. Solving the above linear equations, we get the expressions in \eqref{eq:case3}. Given that $G(x)$ is full row rank, thus $\tilde G^+(x)$ is also full row rank, which leads to the existence of $(A_{\mathcal{I}}\tilde G^+(x){A_{\mathcal{I}}}^\top)^{-1}$ provided that $A_{\mathcal{I}}$ is full row rank.

Clearly, the aforementioned three cases cover all the possible situations when solving problem \eqref{eq:mpop}. The following theorem states the Lipschitz continuous property of the optimal safe control input $u^*(x)$.

\begin{theorem}\label{th:conlipschitz}
Consider problem \eqref{eq:mpop}. Suppose that this problem is feasible for any $x\in\mathcal{C}$. Then the optimal control law $u^*(x)$ is locally Lipschitz continuous in $\mathcal{C}$.
\end{theorem}
\begin{proof}
The multipliers $\lambda(x),\mu(x)$, and control law $u^*(x)$ are bounded because of the relative degree assumption of $B(x)$. Hence, $u^*(x)$ is Lipschitz continuous within every domain in the three cases with expressions \eqref{eq:case1}, \eqref{eq:case2}, \eqref{eq:case3}, since the composition and product of locally Lipschitz continuous functions are Lipshitz continuous. Besides, $u^*(x)$ is well defined and continuous on the boundary of the domain provided that all the constraints satisfy strict complementary slackness and are linear independent \cite{bemporad2002explicit}. Hence, we conclude that the optimal safe controller $u^*(x)$ is locally Lipschitz continuous over $\mathcal{C}$.
\end{proof}
The results in Equations \eqref{eq:case1} -- \eqref{eq:case3} characterize the solutions only in a neighborhood of a specific $x$. Given $x\in\mathcal{C}$, the question remaining is to which case this belongs to. The following section provides an approach to explore the state-space. 

\subsection{Partitioning The State Space}
\label{sec:partition}
After presenting the explicit solution evaluation for all possible cases in the previous section, we now show how to partition the state space and how to address the primal degenerate and dual degenerate cases.

The exploration of the parameter space has been widely investigated in the literature. Our method lies in the class of active-set approaches \cite{feller2013improved}, which tend to enumerate all possible active-inactive combinations of the constraints. The state-space is partitioned into a finite number of critical regions. One difference is that we no longer need to solve additional optimality and feasibility problems as in \cite{pappas2021exact}.

The proposed approach is based on the enumeration of all possible combinations of active and inactive control limit constraints. Let $\mathcal{AS}$ denote the set of active sets, and $\bar{\mathcal{AS}}$ denote the inactive sets
\begin{subequations}\label{eq:enumeration}
\begin{equation}
    \mathcal{AS}:=\{\mathcal{I}_1=\{1\},\ldots,\mathcal{I}_{2^p-1}=\{1,\ldots,p\}\},
\end{equation}
\begin{equation}
    \bar{\mathcal{AS}}:=\{\bar{\mathcal{I}}_1=\{2,\ldots,p\},\ldots,\bar{\mathcal{I}}_{2^p-1}=\emptyset\}.
\end{equation}
\end{subequations}
Our approach is shown in Algorithm \ref{al:partition}. Critical regions for Case 1 and Case 2 are partitioned in lines 1-2, denoted by CR$_1$ and CR$_2$, respectively. Lines 5-6 deal with the linearly dependent constraints. If $A_{\mathcal{I}}$ is not full row rank, this directly indicates that there are linear dependent constraints. This case can be pruned since there are redundant linear equations involved. Line 8 corresponds to Case 3. We note here that in the algorithm, strict inequalities are used to partition the state space, and therefore the boundary of each region is not well defined. Since $u^*(x)$ is continuous on every boundary, one can assign the limit points to any neighbor region . More specifically, the boundary corresponds to those $x$ where the optimisation problem \eqref{eq:mpop} is dual degenerate.
	\begin{algorithm}[htb]
		\small
		\LinesNumbered
		\caption{State space partition algorithm}\label{al:partition}
		\KwIn{Matrix $A,b$, functions $G(x),F(x),\Lambda(x)$, vector $u^{\mathrm{des}}$
		}
		\KwOut{Optimal control law $u^*(x)$, critical regions $\mathrm{CR}_1,\ldots,\mathrm{CR}_l$
		\vspace{3pt}	
		}
		If the control Barrier function constraint is inactive, obtain the critical region by substituting the parametric expressions on the inactive control Barrier function constraint $\mathrm{CR}_1:=\mathcal{C}\cap\{x|F(x)+G(x)u^*(x)+\Lambda(x)<0\}$, where $u^*(x)$ is given by \eqref{eq:case1}\\
		If the control Barrier function constraint is active and the control limits constraints are inactive, the critical region is defined by
		$\mathrm{CR}_2:=\mathcal{C}\cap \{x|Au^*(x)+b<0\}\cap\{x|\lambda(x)>0\}\}$, where $u^*(x)$ and $\lambda(x)$ are given by \eqref{eq:case2}.\\
		If the control Barrier function constraint is active, enumerate all possible pair-wise combinations of active sets $\mathcal{AS}$ and inactive sets $\bar{\mathcal{AS}}$.\\
		\For{i=1:p}
		{
		\eIf{$A_{\mathcal{I}_i}$ is not full row rank}
		{Prune this case}
		{
		The critical region is defined by $\mathrm{CR}_{i+2}:=\mathcal{C}\cap\{x|A_{\bar{\mathcal{I}_i}}u^*(x)+b<0\}\cap\{x|\lambda(x)>0\}\cap\{x|\mu(x)>0\}$, where
		$u^*(x)$, $\lambda(x)$, $\mu(x)$  are given by \eqref{eq:case3}.
		}
		}
	\end{algorithm}

\subsection{Adaptive Control Barrier Function-Based QP}
In this section, we consider the feasibility problem discussed in the previous section. Lipschitz continuity of the controller is guaranteed by the relaxation term, i.e. the class-$\mathcal{K}$ function $\alpha(\cdot)$ in \eqref{eq:cbf}. At some $x\in\mathcal{C}$, the problem could be infeasible due to the control limits or unsuitable $\alpha(\cdot)$. Numerous adaptive control Barrier function approaches have been proposed to tune $\alpha(\cdot)$ dynamically to improve feasibility \cite{xiao2021adaptive}. However none of these results reveals the explicit adapted relaxation term design and demonstrates the Lipschitz continuity of $u^*(x)$. In this section we address these problems by solving an adaptive control Barrier function-based QP problem explicitly.

The adaptive control Barrier function-based QP can be formulated as
\begin{equation}\label{eq:adcbf}
    \begin{split}
    (s^*(x),u^*(x))&=\arg\min_{s,u}~\frac{p_s}{2}(s-1)^2+\frac{1}{2}||u-u^{\mathrm{des}}||_2^2\\
    &\mathrm{s.t.}~F(x)+G(x)u+s\Lambda(x)\le 0,\\
    &~~~~~Au+b\le 0,
    \end{split}
\end{equation}
where $s$ is an adapted parameter to tune the relaxation term dynamically, and $p_s\in\mathbb{R}^+$ is a predefined penalty coefficient. The Lagrangian of the problem \eqref{eq:adcbf} is
\begin{equation}\label{eq:adcbflagrangian}
\begin{split}
 L(s,u,\lambda,\mu) = &\frac{p_s}{2}(s-1)^2+\frac{1}{2}||u-u^{\mathrm{des}}||_2^2+\\\mu^\top (Au+b)
 +&\lambda^\top(F(x)+G(x)u+s\Lambda(x)).
\end{split}
\end{equation}
For every $x$ the problem is a convex QP over both $s$ and $u$, strong duality holds. The KKT conditions for the optimisation problem are 
\begin{subequations}\label{eq:adcbfkkt}
\begin{align}
    &p_ss^*(x)-p_s+\Lambda(x)\lambda=0,\label{eq:adcbfkkt1}\\
    &A^\top \mu+G(x)^\top\lambda+u^*(x)-u^{\mathrm{des}} = 0,\label{eq:adcbfkkt2}\\
    &\mu^\top(Au^*(x)+b)=0,\label{eq:adcbfkkt3}\\
    &\lambda(F(x)+G(x)u^*(x)+s^*(x)\Lambda(x))=0,\label{eq:adcbfkkt4}\\
    &\lambda\ge0,F(x)+G(x)u^*(x)+s^*(x)\Lambda(x)\le0,\label{eq:adcbfkkt4}\\
    &\mu\ge0,Au^*(x)+b\le0.\label{eq:adcbfkkt6}
\end{align}
\end{subequations}

Following the steps in Section \ref{sec:cbfqp}, the problem can be solved by considering three cases. One difference here is that the feasibility requirement is no longer assumed when solving the problem explicitly.
\begin{theorem}\label{th:feasibility}
The optimisation problem \eqref{eq:adcbf} is feasible for any $x\in\mathcal{C}$, and any class-$\mathcal{K}$ function $\Lambda(x)$ if $B(x)$ is a control Barrier function.
\end{theorem}

\begin{proof}
We first prove the theorem for $\Lambda(x)>0$ and then $\Lambda(x)=0$. For any $x'\in\{x|\Lambda(x)>0\}\subseteq \mathcal{C}$, $u'\in\mathcal{U}$ the search space $s\le-\frac{F(x)+G(x)u'}{\Lambda(x)}$ is nonempty since $-\frac{F(x)+G(x)u'}{\Lambda(x)}$ is a finite real scalar. Thus, problem \eqref{eq:adcbf} is feasible. Consider now  the case where $x'\in\{x|\Lambda(x)=0\}\subseteq \mathcal{C}$. According to the definition of control Barrier functions, we have $\sup_{u\in\mathcal{U}}[\mathcal{L}_fB(x)+\mathcal{L}_gB(x)u+1\cdot\Lambda(x)]|_{x=x'}$, the problem \eqref{eq:adcbf} is feasible with $s^*(x)=1$. Therefore, we conclude that for any $x\in\mathcal{C}$, and any class-$\mathcal{K}$ function $\Lambda(x)$, the adaptive control Barrier function based QP problem \eqref{eq:adcbf} is feasible.
\end{proof}

\begin{remark}\label{rem:infeasible}
For the case where ${\mathcal{B}}\subseteq \mathcal{S}$ but $B(x)$ is not a proper Control Barrier Function, \eqref{eq:adcbf} is feasible for any $x\in \mathrm{Int}(\mathcal{B})$. For somer $x\to\partial{\mathcal{B}}$, we directly have that $s^*(x)\to\infty$ if $F(x)+G(x)u<0$ for any $u\in\mathcal{U}$, since $\Lambda(x)\to 0$.
\end{remark}

\textit{Case 1}: Given $x$, the control Barrier function constraint is inactive i.e. $F(x)+G(x)u^*(x)+s\Lambda(x)<0$ and $\lambda=0$. The minimized adapted parameter $s^*(x)$ of \eqref{eq:adcbf} in this case is given by
\begin{equation}\label{eq:adcase1}
    s^*(x)=1.
\end{equation}
The result follows directly from \eqref{eq:adcbfkkt1} by substituting $\lambda =0$.

\textit{Case 2}: Given $x$, the control Barrier function constraint is active, and the control limit constraints are inactive, i.e. $Au+b=0$ and $\mu=0$.  Then the optimal control law $u^*(x)$, the minimized adapted parameter $s^*(x)$, and the Lagrange multiplier $\lambda(x)$ is given by 
\begin{subequations}\label{eq:adcase2}
\begin{align}
    u^*(x)&=u^{\mathrm{des}}-G(x)^\top\lambda(x),\label{eq:adcase2eq1}\\
    s^*(x)&=1-\frac{\Lambda(x)}{p_s}\lambda(x),\label{eq:adcase2eq2}\\
    \lambda(x)&=\frac{p_s(\Lambda(x)+G(x)u^{\mathrm{des}}+F(x))}{p_sG(x) G(x)^\top+\Lambda(x)^2}.\label{eq:adcase2eq3}
\end{align}
\end{subequations}
To prove this, let $\mu=0$ into \eqref{eq:adcbfkkt2}, we obtain the following linear equations
\begin{subequations}\label{eq:adcbflinear2}
\begin{align}
    &p_ss^*(x)-p_s+\Lambda(x)\lambda(x)=0,\\
    &G(x)^\top\lambda(x)+u^*(x)-u^{\mathrm{des}}=0,\\
    &F(x)+G(x)u^*(x)+s^*(x)=0.
\end{align}
\end{subequations}
Solving these linear equations leads to \eqref{eq:adcase2}.

From the results for Case 2 we can see that the adapted parameter $s^*(x)=1$ only when $\Lambda(x)+G(x)u^{\mathrm{des}}+F(x)=0$, which shows that $\lambda(x)=0$. Therefore, we conclude that the adapted parameter needs to be tuned if and only if the control Barrier function constraint is active.

\textit{Case 3}: Given $x$, the control Barrier function constraint is active and the control limit constraints are active with index set $\mathcal{I}$. Then the optimal adapted parameter $s^*(x)$, optimal control input $u^*(x)$ and the Lagrange multipliers $\lambda(x),\mu(x)$ are given by
    \begin{subequations}\label{eq:adcase3}
    \begin{align}
        s^*(x)&=1-\frac{\Lambda(x)}{p_s}\lambda(x),\label{eq:adcase3s}\\
        u^*(x)&=u^{\mathrm{des}}-{A_{\mathcal{I}}}^\top \mu(x)-G(x)^\top \lambda(x),\label{eq:adcase3u}\\
        \mu(x)&= (A_{\mathcal{I}}{A_{\mathcal{I}}}^\top)^{-1}(A_{\mathcal{I}}u^{\mathrm{des}}-A_{\mathcal{I}}G(x)^\top\lambda(x)+b),\label{eq:adcase3mu}\\
        \lambda(x)&= (p_sG(x)\tilde A^-_{\mathcal{I}}G(x)^\top+\Lambda(x)^2)^{-1}\label{eq:adcase3lambda}\\
        &\times p_s(F(x)+G(x)\tilde A^-_{\mathcal{I}}u^{\mathrm{des}}+\Lambda(x)\nonumber\\
        &-G(x){A_{\mathcal{I}}}^\top(A_{\mathcal{I}}{A_{\mathcal{I}}}^\top)^{-1}b_{\mathcal{I}}),\nonumber
    \end{align}
    \end{subequations}
where $\tilde A^-_{\mathcal{I}} := I-{A_{\mathcal{I}}}^\top(A_{\mathcal{I}}{A_{\mathcal{I}}}^\top)^{-1}A_{\mathcal{I}}$,
The results are obtained by solving the following linear equations
\begin{subequations}
\begin{align}\label{eq:adlinearkkt}
     &p_ss^*(x)-p_s+\Lambda (x)\lambda(x) =0,\\
    &{A_{\mathcal{I}}}^\top \mu(x)+G(x)^\top\lambda(x)+u-u^{\mathrm{des}}=0,\\
    &F(x)+G(x)u+s^*(x)\Lambda(x)=0,\\
    &A_{\mathcal{I}}u+b_{\mathcal{I}}=0,
\end{align}
\end{subequations}
with the observation that $p_sG(x)\tilde A^-_{\mathcal{I}}G(x)^\top+\Lambda(x)^2$ is non-zero, since $\tilde A^-_{\mathcal{I}}\succeq 0$ and $p_s>0$, $\Lambda(x)\ge 0$.

Having listed all the possible cases for explicit adapted parameter design, the region characteristic for this problem follows similar steps as in Algorithm \ref{al:partition}. Ideally, $s^*(x)=1$ in the feasible critical regions for the original problem \eqref{eq:mpop}. However, this cannot be achieved unless we have  $p_s\to\infty$.

The following theorem states the Lipschitz continuity proof for controller design with adapted formulation\ref{th:adconlipschitz}.

\begin{theorem}\label{th:adconlipschitz}
The optimal control law $u^*(x)$ obtained by solving \eqref{eq:adcbf} is locally Lipschitz continuous in $\mathcal{C}$. 
\end{theorem}
\begin{proof}
The proof is similar to that of Theorem \ref{th:conlipschitz}, since $s^*(x)$ is Lipschitz continuous in every critical region, and is continuous on the boundary.
\end{proof}
The new relaxation term $s^*(x)\alpha(B(x))$ has no obvious monotonicity property according to \eqref{eq:adcase2eq2} and \eqref{eq:adcase3s}. In fact any Lipschitz continuous function which is zero for any $x\in\{x|B(x)=0\}$ guarantees safety as well as renders the controller Lipschitz continuous.

\section{Simulation Results}\label{sec:simulation}
We first illustrate the state space exploration for explicit safe controller synthesis. Consider the following linear system
\begin{equation}\label{eq:linexample}
    \begin{bmatrix}
    \dot x_1\\
    \dot x_2
    \end{bmatrix}
    =
    \begin{bmatrix}
    1&2\\
    1&1
    \end{bmatrix}
    \begin{bmatrix}
    x_1\\
    x_2
    \end{bmatrix}
    +
    \begin{bmatrix}
    u_1\\
    u_2
    \end{bmatrix},
\end{equation}
with affine control limits $-1\le u_1\le 1$, $-1\le u_2\le 1$. A quadratic control Barrier function is defined by $B(x)=-x_1^2-x_2^2+9$. The class-$\mathcal{K}$ function is $\alpha(B(x))=0.5B(x)$, and $u^{\mathrm{des}}=[0.5,0.5]^\top$. Following the state-space partitioning procedure in Algorithm \ref{al:partition}, we obtain the critical regions listed in Table \ref{tab:cbf}. We note here that the active index set $\mathcal{I}$ contains at most one element in this case due to the linear independence requirement. 

\begin{figure*}[t]
    \centering
    \subfigure[domain partition for $u^*(x)$ derived from \eqref{eq:cbf}]{
        \includegraphics[width=0.3\textwidth]{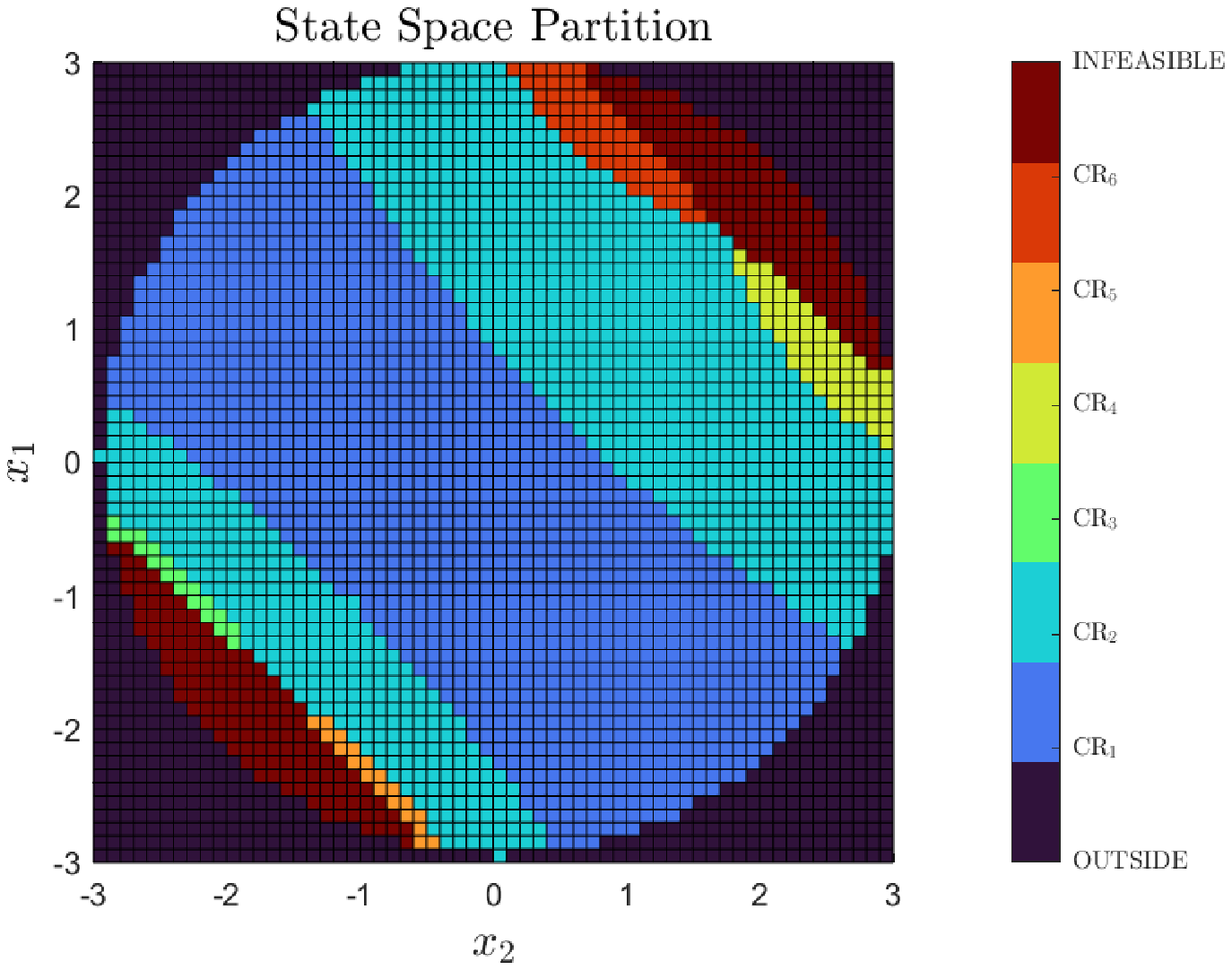}
        \label{fig:space1}
    }
    \subfigure[value of $u_1^*(x)$]{
        \includegraphics[width=0.3\textwidth]{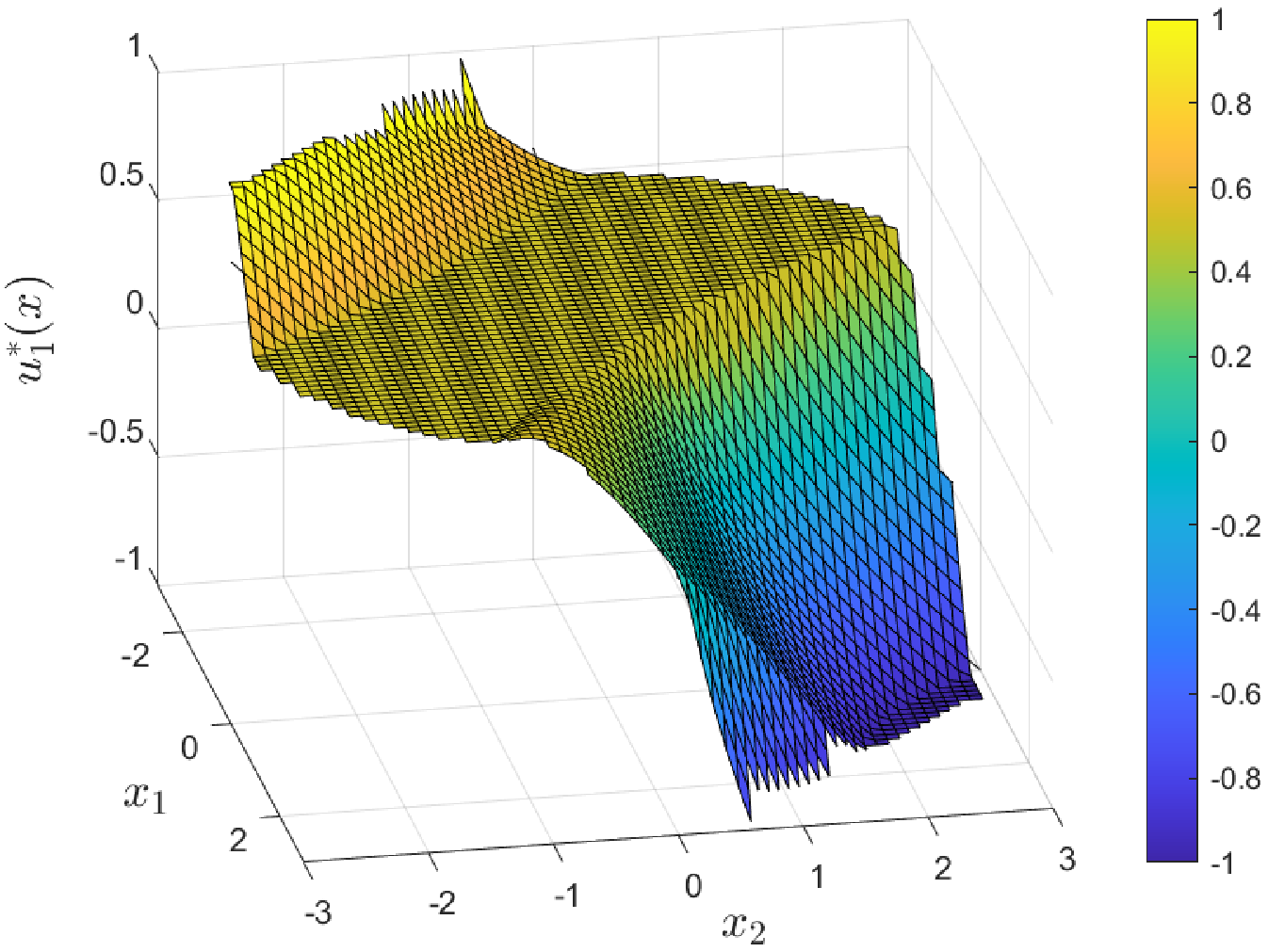}
        \label{fig:space1u1}
    }
    \subfigure[value of $u_2^*(x)$]{
        \includegraphics[width=0.3\textwidth]{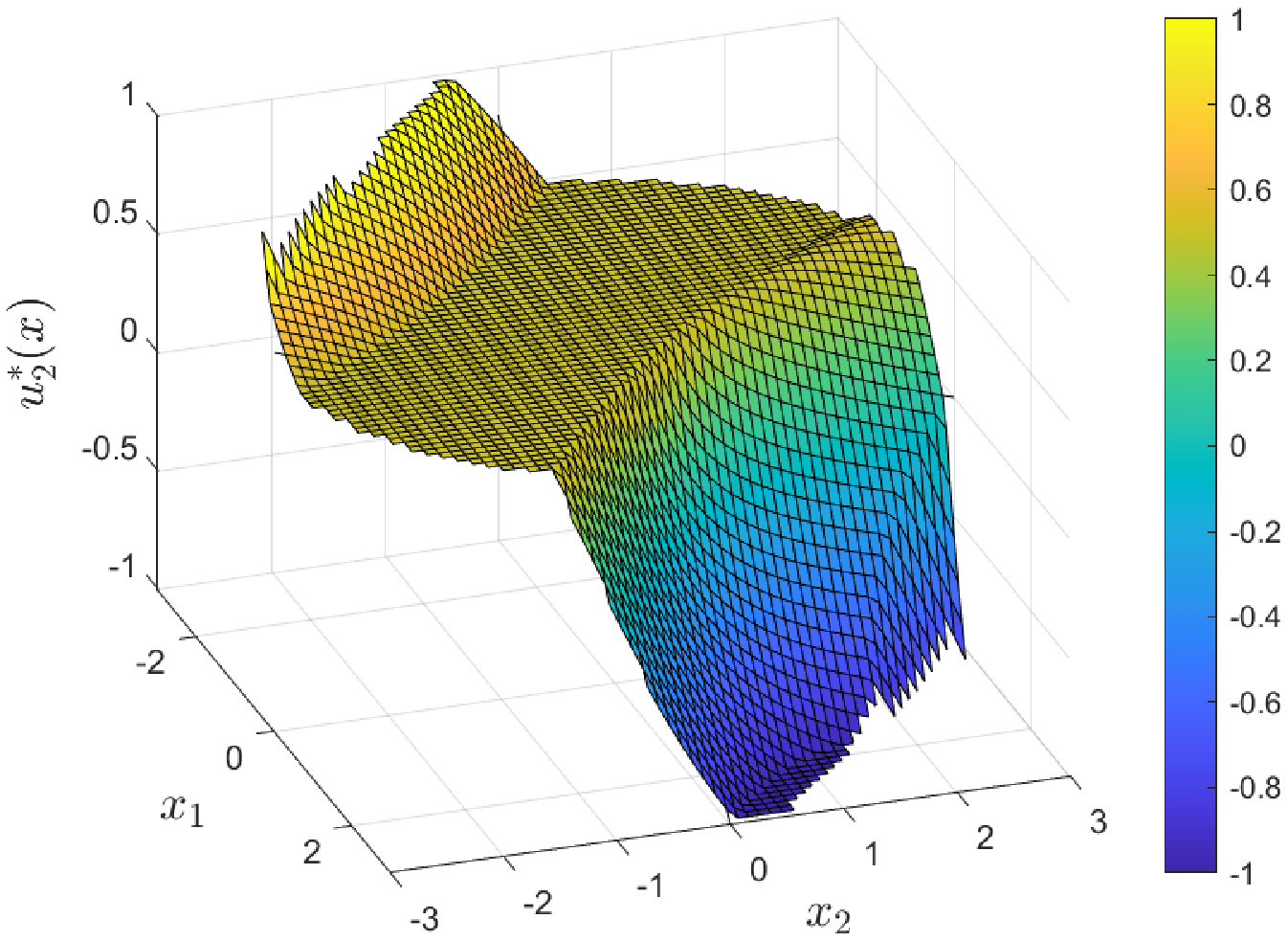}
        \label{fig:space1u2}
    }
    \caption{Value of $u^*(x)$ in the state space, the piece-wise explicit expression is Lipschitz continuous}
    \label{fig:cbf}
\end{figure*}

\begin{figure*}[t]
    \centering
    \subfigure[domain partition for $u^*(x)$ derived from \eqref{eq:cbf}]{
        \includegraphics[width=0.23\textwidth]{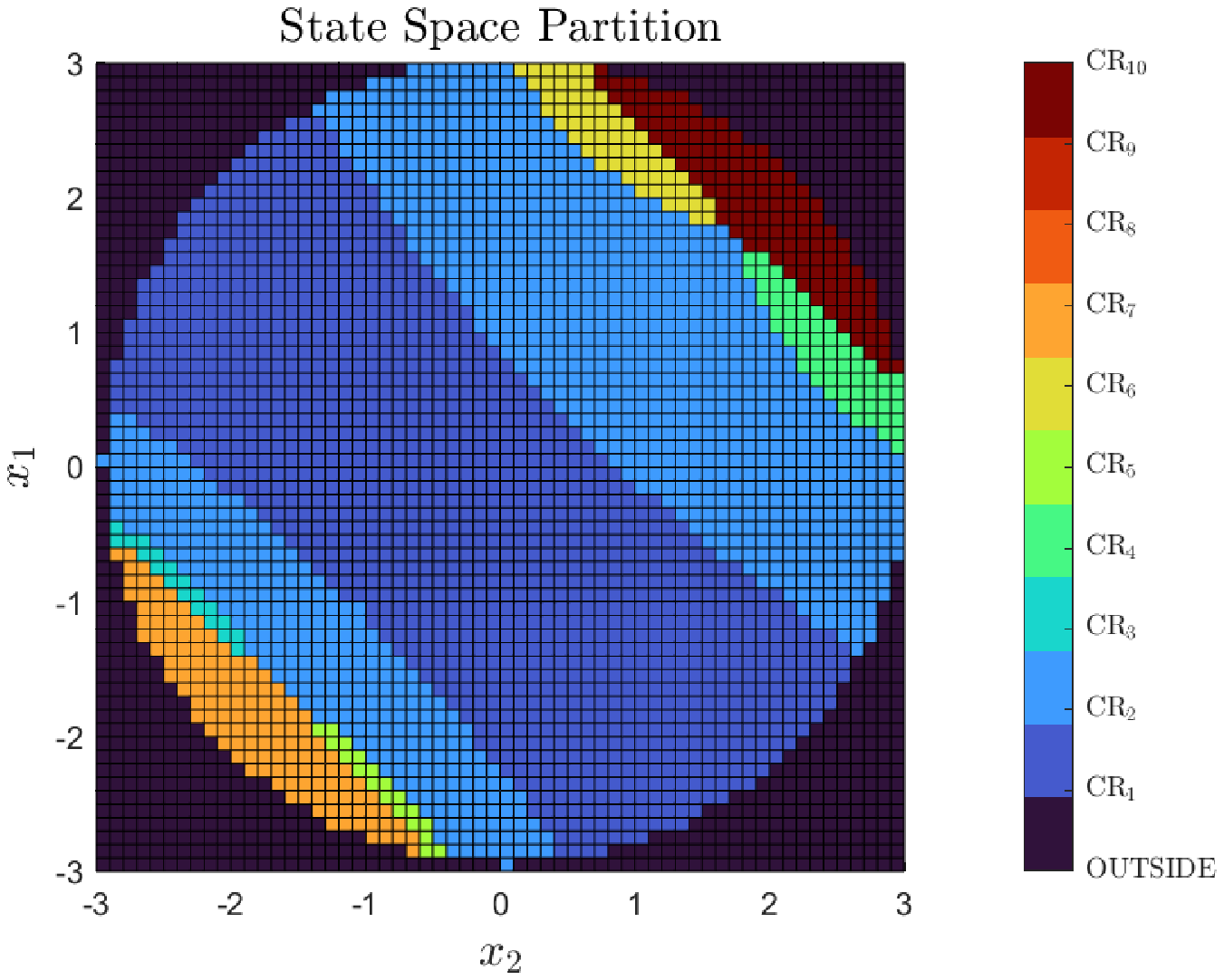}
        \label{fig:adspace2}
    }
    \subfigure[value of $u_1^*(x)$]{
        \includegraphics[width=0.23\textwidth]{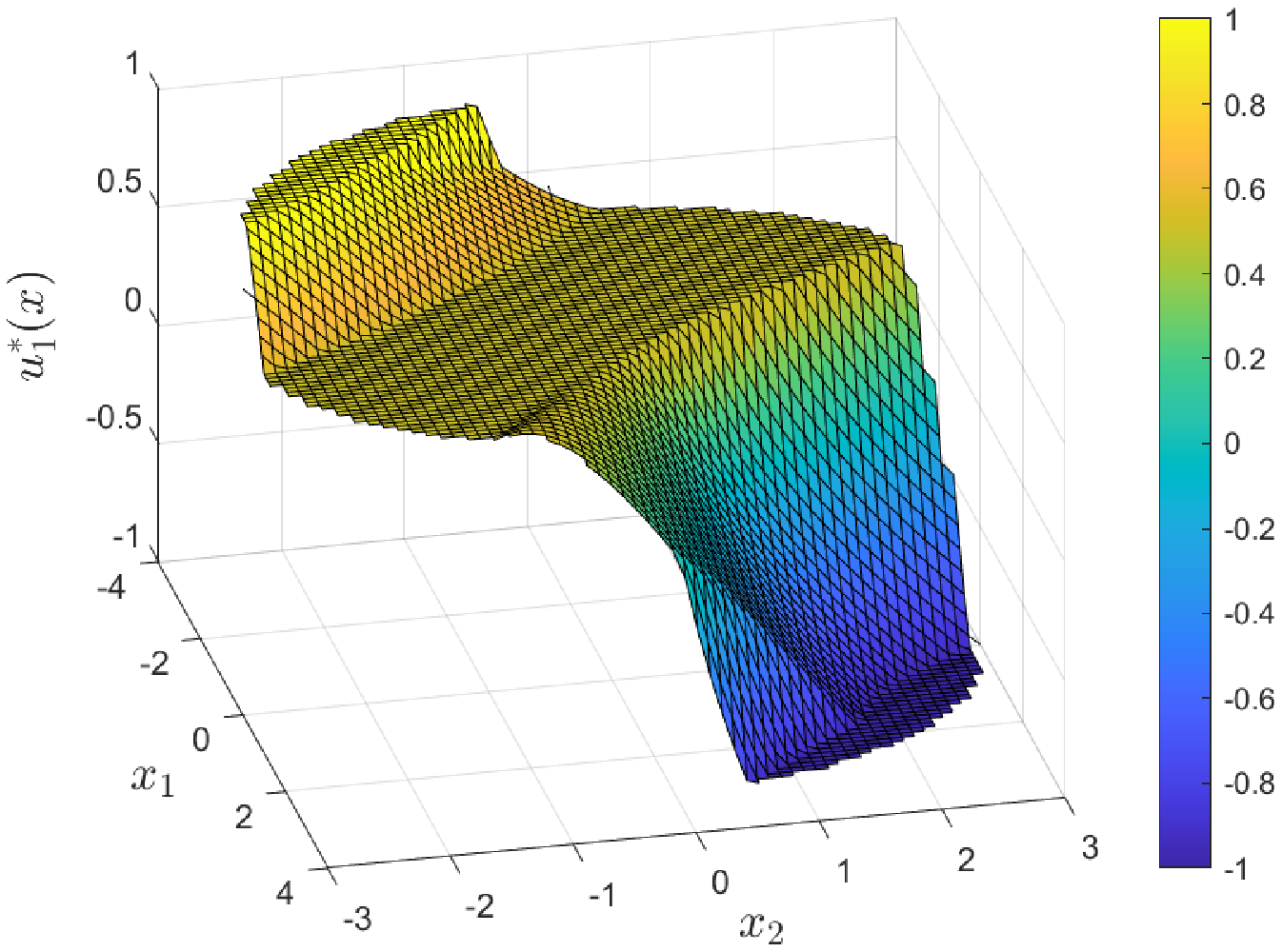}
        \label{fig:adspace2u1}
    }
    \subfigure[value of $u_2^*(x)$]{
        \includegraphics[width=0.23\textwidth]{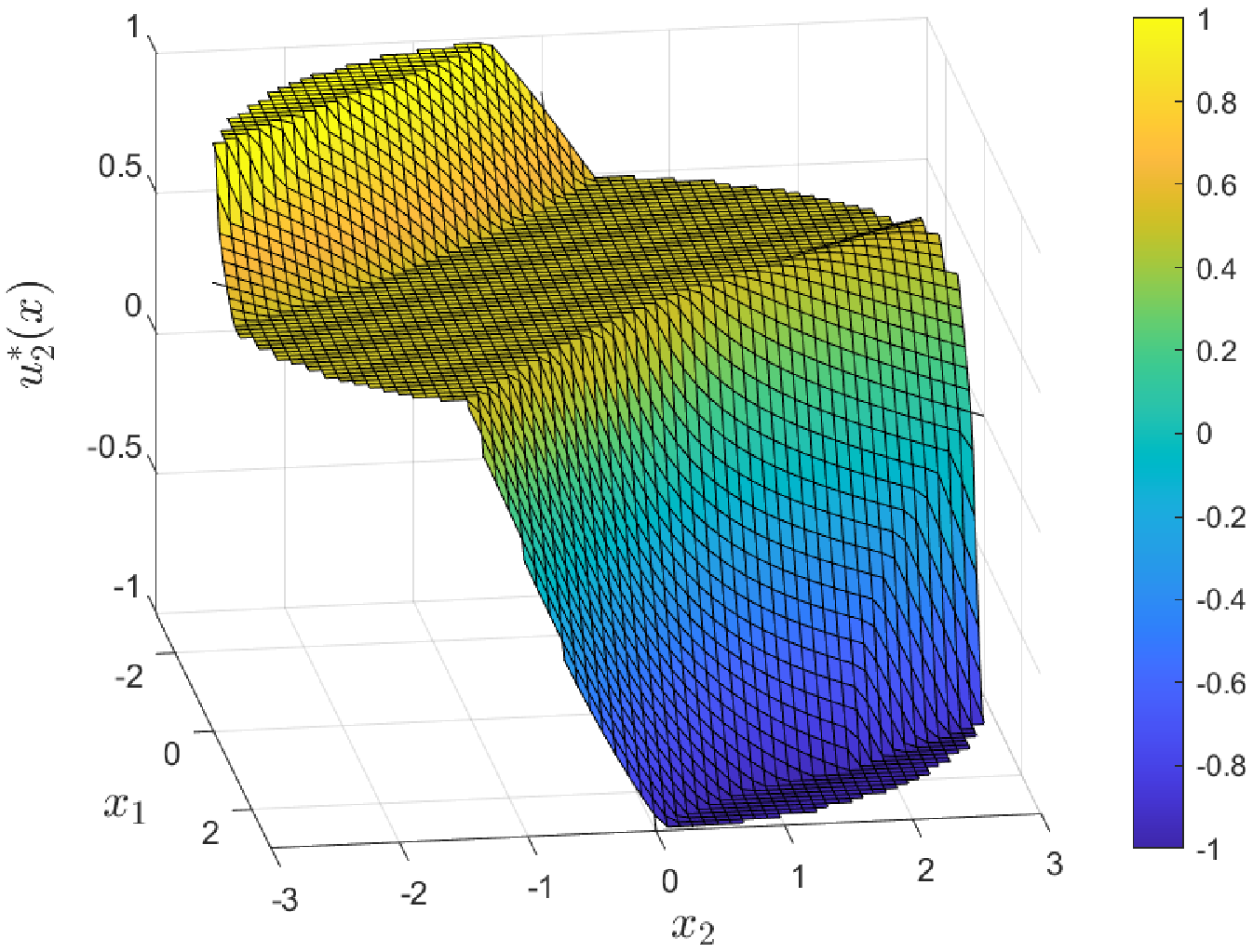}
        \label{fig:adspace2u2}
    }
    \subfigure[value of $s^*(x)$]{
        \includegraphics[width=0.23\textwidth]{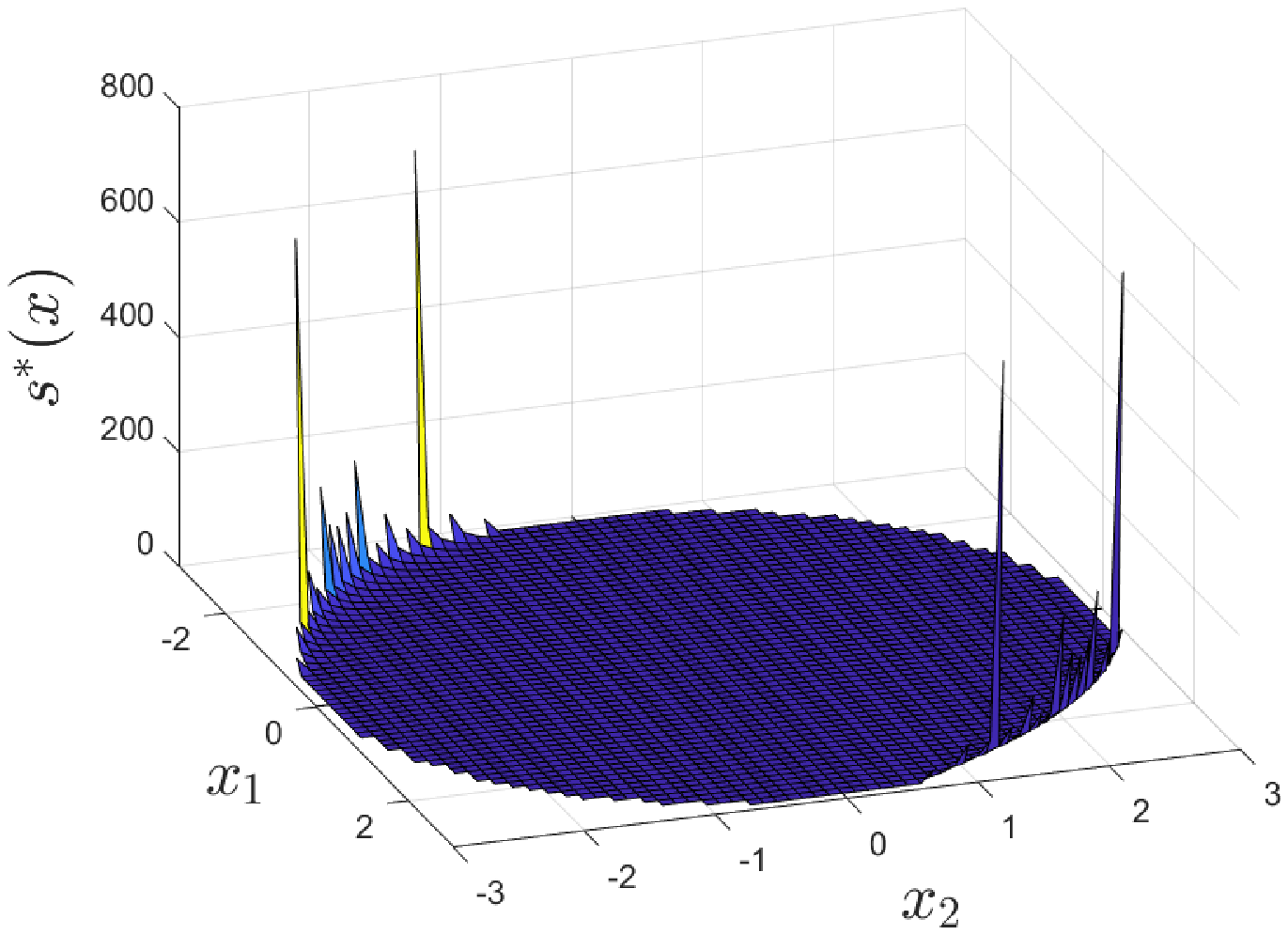}
        \label{fig:adspace2s}
    }
    \caption{Value of $u^*(x)$ in the state space, the piece-wise explicit expression is Lipschitz continuous}
    \label{fig:adcbf}
\end{figure*}

\begin{table}[h]
\centering
\begin{threeparttable}
\begin{tabular}{ccc}
\toprule[1.5pt]
Critical Regions & \makecell{Conditions Activeness} & Index Set $\mathcal{I}$\\
\midrule[0.75pt]
\makecell{CR1} & \makecell{$F(x)+G(x)u+\Lambda(x)<0$\tnote{*}} & \makecell{$\emptyset$}\\
\hline
\makecell{CR2} & \makecell{$F(x)+G(x)u+\Lambda(x)=0$}&\makecell{$\emptyset$}\\ \hline
\makecell{CR3} & \makecell{$F(x)+G(x)u+\Lambda(x)=0$\\$u_1(x)=1$}&[1]\\ \hline
\makecell{CR4} & \makecell{$F(x)+G(x)u+\Lambda(x)=0$\\$u_1(x)=-1$}&[2]\\ \hline
\makecell{CR5} & \makecell{$F(x)+G(x)u+\Lambda(x)=0$\\$u_2(x)=1$}&[3]\\ \hline
\makecell{CR6} & \makecell{$F(x)+G(x)u+\Lambda(x)=0$\\$u_2=-1$}&[4]\\
\hline
\makecell{CR7} &
\makecell{NaN\tnote{**}}&NaN\\
\bottomrule[1.5pt]
\end{tabular}
\begin{tablenotes}
        \item[*] The control Barrier function constraint is inactive (Case 1)
        \item[**] The problem is infeasible
      \end{tablenotes}  
\end{threeparttable}

\caption{State space partitioning for the control Barrier functions based quadratic programming}
\label{tab:cbf}
\end{table}

Figure \ref{fig:space1} shows the partitioned regions. The safe region is partitioned into seven critical regions, in which $u^*(x)$ is defined as a piece-wise continuous function. The problem is infeasible inside some regions of the safe region which suggests that either $B(x)$ is not a candidate control Barrier function, or the $\alpha(B(x))=0.5B(x)$ is not a proper relaxation function in class-$\mathcal{K}$. Figures \ref{fig:space1u1}-\ref{fig:space1u2} show the value of $u_1^*(x)$ and $u_2^*(x)$ in the state space. In CR$_1$, $u^*(x)=u^{\mathrm{des}}$ has a constant value. Overall $u^*(x)$ is locally Lipschitz continuous. 
\begin{table}[h]
\centering
\begin{threeparttable}
\begin{tabular}{ccc}
\toprule[1.5pt]
Critical Regions & \makecell{Conditions Activeness} & Index Set $\mathcal{I}$\\
\midrule[0.75pt]
\makecell{CR1} & \makecell{$F(x)+G(x)u+s\Lambda(x)<0$} & \makecell{$\emptyset$}\\
\hline
\makecell{CR2} & \makecell{$F(x)+G(x)u+s\Lambda(x)=0$}&\makecell{$\emptyset$}\\ \hline
\makecell{CR3} & \makecell{$F(x)+G(x)u+s\Lambda(x)=0$\\$u_1(x)=1$}&[1]\\ \hline
\makecell{CR4} & \makecell{$F(x)+G(x)u+s\Lambda(x)=0$\\$u_1(x)=-1$}&[2]\\ \hline
\makecell{CR5} & \makecell{$F(x)+G(x)u+s\Lambda(x)=0$\\$u_2(x)=1$}&[3]\\ \hline
\makecell{CR6} & \makecell{$F(x)+G(x)u+s\Lambda(x)=0$\\$u_2=-1$}&[4]\\
\hline
\makecell{CR7} &
\makecell{$F(x)+G(x)u+s\Lambda(x)=0$\\$u_1=1$\\$u_2=1$}&[1,3]\\
\hline
\makecell{CR8} &
\makecell{$F(x)+G(x)u+s\Lambda(x)=0$\\$u_1=1$\\$u_2=-1$}&[1,4]\\
\hline
\makecell{CR9} &
\makecell{$F(x)+G(x)u+s\Lambda(x)=0$\\$u_1=-1$\\$u_2=1$}&[2,3]\\
\hline
\makecell{CR10} &
\makecell{$F(x)+G(x)u+s\Lambda(x)=0$\\$u_1=-1$\\$u_2=-1$}&[2,4]\\
\bottomrule[1.5pt]
\end{tabular}
\end{threeparttable}

\caption{State space partitioning for the adaptive control Barrier functions based quadratic programming}
\label{tab:adcbf}
\end{table}

As this problem is infeasible close to the boundary of the safe region, we can use the adaptive control Barrier functions formulation \eqref{eq:adcbf} to improve feasibility. With a larger coefficient, the control Barrier function constraint is relaxed and the problem is rendered solvable if $B(x)$ is a candidate Control Barrier Function according to Theorem \ref{th:feasibility}. Following the same procedure as in Algorithm \ref{al:partition}, the state space partitioning is shown in Table \ref{tab:adcbf}. Figure \ref{fig:adcbf} shows the domain partition, value of $u^*(x)$ and $s^*(x)$. It can be seen that the problem is feasible for any $x\in\partial {\mathcal B}$ in Figure \ref{fig:adspace2}. Figures \ref{fig:adspace2u1}-\ref{fig:adspace2u2} show that $u^*(x)$ is locally Lipschitz continuous with the adaptive Control Barrier Function formulation. In Figure \ref{fig:adspace2s}, note that $s^*(x)\to\infty$ when $x\to\partial \mathcal{B}$. This indicates that $B(x)$ is not a candidate Control Barrier Function for system \eqref{eq:linexample} with control limits.

\section{Conclusion}\label{sec:conclusion}
In this paper we investigated the explicit safe controller synthesis problem. The proposed approach was based on parameterized control Barrier functions-based quadratic programming. Exploring the state space, multiple disjoint critical regions are identified as domains for piece-wise explicit controller design. For the case where the problem is infeasible we propose an adaptive coefficient adaptation scheme. Simulation results demonstrate our results. In the future we will explore how to enhance the smoothness of the resulting controller with carefully designed relaxation terms. 

\bibliographystyle{ieeetr}
\bibliography{ref.bib}
\end{document}